\documentclass[english,11pt]{amsart}
\usepackage{amsmath}
\usepackage{amsthm}
\usepackage{pdfsync}
\usepackage{a4wide}
\usepackage{hyperref}

\usepackage{amssymb}
\usepackage{amsfonts}
\usepackage{color}
\usepackage{enumerate}
\usepackage{graphicx}
\usepackage{amsfonts}
\usepackage{color}
\usepackage{enumerate}

\newtheorem{theorem}{Theorem}

\newtheorem{remark}{Remark}

\def\rr{\mathbb{R}}
\def\intr{\int_{\rr^N}}

\def\n{\nabla}
\def\nn{\mathbb{N}}
\def\eps{\epsilon}

\begin{document}

\title[Short proof of the Hardy-Rellich inequality]{A new proof of the Hardy-Rellich inequality in any dimension}
\author[C. Cazacu]{Cristian Cazacu}
\address[C.  Cazacu]{\\Faculty of Mathematics and Computer Science \& ICUB,  University of Bucharest  \\
14 Academiei Street \\ 010014 Bucharest\\ Romania\\
}
\email{cristian.cazacu@fmi.unibuc.ro}

\begin{abstract}
The Hardy-Rellich inequality in the whole space with the best constant was firstly proved by Tertikas and Zographopoulos in Adv. Math. (2007) in higher dimensions $N\geq 5$. Then it was extended to lower dimensions $N\in \{3, 4\}$ by Beckner in Forum Math. (2008) and Ghoussoub-Moradifam in Math. Ann. (2011) by applying totally different techniques. 

  In this note we refine the method implemented by Tertikas and Zographopoulos, based on spherical harmonics decomposition, to give an easy and compact proof of the optimal Hardy-Rellich inequality in any dimension $N\geq 3$.   In addition, we provide minimizing sequences which were not explicitly mentioned in the quoted papers, emphasizing their symmetry breaking in lower dimensions $N\in \{3,4\}$. We also show that the best constant is not attained in the proper functional space.
\end{abstract}

\maketitle



In this note we first present a new  unified proof for the following well-known optimal Hardy-Rellich  inequality in any dimension.
\begin{theorem}\label{HardyBilaplacian} Assume $N\geq 3$. Then, for any $u\in C_c^\infty(\rr^N)$ it holds
\begin{equation}\label{optconst}
\intr |\Delta u|^2 dx \geq C(N) \intr \frac{|\n u|^2}{|x|^2} dx,
\end{equation}
where
\begin{equation}
C(N):=\left\{\begin{array}{ll}
              \frac{N^2}{4}, & N\geq 5 \\[5pt]
              3, & N=4 \\[5pt]
              \frac{25}{36}, & N=3. \\[5pt]
            \end{array}\right.
\end{equation}
\end{theorem}
To the best of our knowledge inequality \eqref{optconst} was firstly analyzed and proved by Tertikas-Zographopoulos \cite{TertiZagra} in higher dimensions $N\geq 5$. Their method applies spherical harmonics decomposition but their proof fails for lower dimensions $N\in \{3, 4\}$.  Soon after that, inequality \eqref{optconst} was firstly   completed in any dimensions $N\geq 3$  by Beckner \cite{MR2431497}, making usage of Fourier transform tools.  Subsequently,  Moradifam-Ghoussoub \cite{MR2753796} developed a quite general theory which allowed them to obtain the most classical functional inequalities and their improvements in the literature. The authors in \cite{MR2753796} combine the method in \cite{TertiZagra} with some ideas from \cite{agmon, MR0374628, MR0342829} reducing the problem to determine positive solutions for some parametric ordinary differential equations of Bessel-type. In particular, the authors in \cite{MR2753796} justify Theorem \ref{HardyBilaplacian}. However, their proof requires to split the analysis into several parts in which they distinguish different techniques in the cases $N\geq 5$ than for $N\in \{3, 4\}$.

We point out that the authors in \cite{MR2753796} considered inequalities in bounded domains but they can be trivially extended to the whole space. It is classical for functional inequalities that the advantage of working in bounded domains allows to improve them by adding positive lower order reminder terms.   It is also worth mentioning the preprint \cite{Futoshi} which complements the above papers with Rellich-type inequalities for vector fields.

The first novelty of this note regards a short (but detailed) and compact proof of Theorem \ref{HardyBilaplacian} in any dimension $N\geq 3$ by means of the spherical harmonics decomposition. In fact, we show that the same  technique applied in  \cite{TertiZagra} to prove Theorem \ref{HardyBilaplacian} for higher dimensions  $N\geq 5$ (but slightly modified computations) could be easily extended to any dimension $N\geq 3$.

Moreover, although the constant $C(N)$ in Theorem \ref{HardyBilaplacian} is optimal, that is
\begin{equation}\label{optim}
  C(N)=\inf_{u\in C_{c}^{\infty}(\rr^N)\setminus \{0\}}\frac{\int_{\rr^N}|\Delta u|^2 dx}{\int_{\rr^N} |\nabla u|^2 /|x|^2dx}, 
\end{equation}
  it seems that the authors in \cite{MR2753796, MR2431497} do not explicitly give minimizing sequences in lower dimensions $N\in \{3, 4\}$ for $C(N)$, see, e.g. \cite[Th. 3.5]{MR2753796} and its proof. However, in \cite{TertiZagra} minimizing sequences are given in dimensions $N\geq 5$.

Next we provide minimizing sequences in the cases $N\in \{3, 4\}$. We also show the non-attainability (in the largest possible Hilbert space) of the best constant $C(N)$ for any $N\geq 3$, fact which was not emphasized in the quoted papers. 

In order to state our results we need some preliminary facts. First let us consider the Hilbert space $\mathcal{D}^{2,2}(\rr^N)$ to be the completion of $C_c^{\infty}(\rr^N)$ in the norm
$$\|u\|=\left(\int_{\rr^N} |\Delta u|^2 dx\right)^{1/2}. $$
Of course, $\|\cdot\|$ is a norm on $C_c^\infty(\rr^N)$ due to the weak maximum principle for harmonic functions. 

 In view of that, the optimization problem \eqref{optim} transfers to the larger space   $\mathcal{D}^{2,2}(\rr^N)$, i.e.
$$  C(N)=\inf_{u\in \mathcal{D}^{2,2}(\rr^N)\setminus \{0\}}\frac{\int_{\rr^N}|\Delta u|^2 dx}{\int_{\rr^N} |\nabla u|^2 /|x|^2dx},$$
 which is the natural space where to look for minimizers.
In addition, we consider a smooth cut-off function $g\in C_c^\infty(\rr)$  such
 $$g(r)=\left\{\begin{array}{ll}
          1, & \textrm{ if } |r|\leq 1 \\
          0, & \textrm{ if } |r|\geq 2.
        \end{array}\right. $$

We claim
\begin{theorem}[Minimizing sequences]\label{minimization}
Let $\eps>0$ and define de sequence \begin{equation}\label{minseq}
            u_\eps(x)=\left\{
            \begin{array}{cc}
              |x|^{-\frac{N-4}{2}+\eps}g(|x|), & \textrm{ if } N\geq 5 \\
              |x|^{-\frac{N-4}{2}+\eps}g(|x|)\phi_1\left(\frac{x}{|x|}
              \right),  & \textrm{ if } N\in \{3, 4\} \\
            \end{array}
            \right.
          \end{equation}
          where $\phi_1$ is a spherical harmonic function of degree 1 such that $\|\phi_1\|_{L^2(S^{N-1})}=1$.
Then $\{u_\eps\}_{\eps>0}\subset \mathcal{D}^{2,2}(\rr^N)$ is a minimizing sequence for $C(N)$, i. e.
\begin{equation}\label{limit-seq}
  \frac{\int_{\rr^N}|\Delta u_\eps|^2 dx}{\int_{\rr^N} |\nabla u_\eps|^2/|x|^2 dx }\searrow C(N), \textrm{ as } \eps \searrow 0.
\end{equation}
Besides,  the constant $C(N)$ is not attained  in $\mathcal{D}^{2,2}(\rr^N)$ (there are no minimizers in $\mathcal{D}^{2,2}(\rr^N)$). 
\end{theorem}

\begin{remark}\label{rem1}
The first part of Theorem \ref{minimization} is relevant for $N\in \{3, 4\}$. 	The fact that $\{u_\eps\}_{\eps>0}$ in \eqref{minseq}, when $N\geq 5$,  is a minimizing sequence is void in view of \cite[Th. 6.6]{TertiZagra} by taking $m=k=0$ and $\phi_0(\sigma)=constant$. Our cut-off function is slightly different than the one in \cite{TertiZagra} but this is not an issue.   
\end{remark}

\section*{Proof of Theorem \ref{HardyBilaplacian}}

The proof follows in several steps as follows.

\subsection*{Step I: Spherical coordinates} We appeal to spherical coordinates instead of cartesian coordinates. The coordinates transformation $ x\in \rr^N \mapsto (r, \sigma) \in (0, \infty)\times S^{N-1}$, where $S^{N-1}$ is the $N-1$-dimensional sphere with respect to the Hausdorff measure in $\rr^N$, is very convenient in $\rr^N$ since we can easily expand in Fourier series. Firstly, let us recall that the expression of the Laplace operator in spherical coordinates is given by
\begin{equation}\label{Laplacian}
\Delta=\partial^2_{rr}+\frac{N-1}{r} \partial_r+\frac{1}{r^2} \Delta_{S^{N-1}},
\end{equation}
where $\partial_r$ and $\partial^2_{rr}$ are both partial derivatives of first and second order with respect to the radial component $r$ whereas  $\Delta_{S^{N-1}}$ represents the Laplace-Beltrami operator with respect to the metric tensor  on  $S^{N-1}$.
Without loss of generality, by density arguments, we may assume  $u\in C_0^\infty(\rr^N\setminus\{0\})$.  Next we apply the spherical harmonics decomposition to expand $u$ as
$$u(x)=u(r\sigma)=\sum_{k=0}^{\infty} u_k(r)\phi_k(\sigma),$$
 It is well-known that such series expansion is possible since there exists an orthogonal basis $\{\phi_k\}_{k\geq 0}$ in $L^2(S^{N-1})$ constituted by spherical harmonic functions $\phi_k$  of degree $k$.  Up to a normalization, we may assume that $\{\phi_k\}_k$ is an orthonormal basis in $L^2(S^{N-1})$. Moreover, such $\phi_k$ are smooth eigenfunctions for the Laplace-Beltrami operator $\Delta_{S^{N-1}}$ with the corresponding eigenvalues $c_k=k(k+N-2)$, $k\geq 0$.   To be more precise, we have the following properties
\begin{equation}\label{orthoeigen}
\left\{
\begin{array}{ll}
  -\Delta_{S^{N-1}}\phi_k=c_k \phi_k  \textrm{ on } S^{N-1}, \\[6pt]
 -\int_{S^{N-1}} \Delta_{S^{N-1}}\phi_k \phi_l d\sigma= \int_{S^{N-1}} \n_{S^{N-1}}\phi_k\cdot \n_{S^{N-1}} \phi_l d\sigma\\[3pt]
  =c_k \int_{S^{N-1}} \phi_k \phi_l d\sigma=c_k \delta_{lk},  \quad k, l\in \nn,
\end{array}
\right.
\end{equation}
where $\delta_{lk}$ represents the Kronecker symbol. Next, we will write $u_k^{\prime}$ and $u_{k}^{\prime\prime}$ to express both first and second derivatives of the Fourier coefficients $\{u_k\}_k$.
In view of the well-known relation
$$|\n u|^2=|\partial_r u|^2 +\frac{|\n_{S^{N-1}}u|^2}{r^2}$$
and the co-aria formula, we express both integrals in \eqref{optconst} in terms of the coefficients $\{u_k\}_k$. Applying the properties \eqref{orthoeigen} we successively obtain
\begin{equation}\label{gradform}
\intr \frac{|\n u|^2}{|x|^2} dx = \sum_{k=0}^{\infty}\left( \int_0^\infty r^{N-3} |u_k^{\prime}|^2 dr + c_k \int_0^\infty r^{N-5} u_k^2 dr\right).
\end{equation}
Moreover, in view of \eqref{Laplacian} we can easily get
\begin{equation}\label{Bil}
\intr |\Delta u|^2 dx=\sum_{k=0}^{\infty} \int_0^\infty r^{N-1}\left(|\Delta_r u_k|^2 +\frac{c_k^2}{r^4}u_k^2 -\frac{2c_k}{r^2} u_k \Delta_r u_k  \right) dr
\end{equation}
where $\Delta_r:=\partial^2_{rr}+\frac{N-1}{r} \partial_r$ is the radial part of the Laplacian in \eqref{Laplacian}. Finally, integration by parts in \eqref{Bil} leads to \begin{multline}\label{finalLap}
\intr |\Delta u|^2 dx=\sum_{k=0}^{\infty} \Big(\int_0^\infty r^{N-1}|u_k^{\prime \prime}|^2 dr + (N-1+2c_k) \int_{0}^{\infty} r^{N-3} |u_k^{\prime}|^2 dr\\
  +\left(c_k^2 + 2c_k (N-4)\right)\int_0^\infty r^{N-5} u_k^2 dr \Big).
\end{multline}
In the sequel,  we prove Theorem \ref{HardyBilaplacian} taking advantage of  identities \eqref{gradform} and \eqref{finalLap}.

\subsection*{Step II: Weighted 1-d Hardy inequalities}
Next, we will apply the following weighted Hardy-Rellich type inequalities
\begin{equation}\label{WeightedRellich}
\int_0^\infty r^{N-1} |u_k^{\prime \prime}|^2 dr \geq \frac{(N-2)^2}{4} \int_0^\infty r^{N-3} |u_k^{\prime}|^2 dr, \quad \forall k\geq 0.
\end{equation}
\begin{equation}\label{weightedhardy}
\int_0^\infty r^{N-3} |u_k^{\prime}|^2 dr \geq \frac{(N-4)^2}{4} \int_0^\infty r^{N-5} u_k^2 dr, \quad \forall k\geq 0.
\end{equation}
The proofs of inequalities \eqref{WeightedRellich} and \eqref{weightedhardy} are straightforward and follow in a similar way. Inequality \eqref{WeightedRellich} is nothing else than the classical Hardy inequality for radial functions but it can be proven independently mimicking the proof of \eqref{weightedhardy}. For the sake of clarity let us give a few lines proof of \eqref{weightedhardy}.  Indeed,
\begin{align}
\int_0^\infty r^{N-5} u_k^2 dr &= \frac{1}{N-4} \int_0^\infty \left(r^{N-4}\right)^{\prime} u_k^2 dr = \frac{-2}{N-4}\int_0^\infty r^{N-4} u_k u_k^{\prime} dr \nonumber\\
&\leq \frac{2}{N-4} \left(r^{N-3} |u_k^{\prime}|^2 dr \right)^{1/2}\left(\int_0^\infty r^{N-5} u_k^2 dr \right)^{1/2},
\end{align}
where the last step is just the Cauchy-Schwarz inequality. Comparing the extreme terms above by taking squares we finally obtain \eqref{weightedhardy}.

\subsection*{Step III: End of the proof} We will make usage of Step I and Step II when comparing both integrals in \eqref{optconst}.

First we split the term on the right hand side in \eqref{finalLap} into the sum $I_1+I_2$ where
$$I_1:=\sum_{k=0}^{\infty} \left(\int_0^\infty r^{N-1}|u_k^{\prime \prime}|^2 dr + (N-1) \int_{0}^{\infty} r^{N-3} |u_k^{\prime}|^2 dr\right)$$
denotes the radial part of the expansion in \eqref{finalLap},  whereas
$$I_2:=\sum_{k=0}^{\infty} \left(2c_k \int_{0}^{\infty} r^{N-3} |u_k^{\prime}|^2 dr+ \left(c_k^2 + 2c_k (N-4)\right)\int_0^\infty r^{N-5} u_k^2 dr\right)$$
is its spherical part.

Then, due to \eqref{WeightedRellich} we have
\begin{equation}\label{eq1}
I_1 \geq \frac{N^2}{4} \sum_{k=0}^{\infty} \int_0^\infty r^{N-3}|u_k^{\prime}|^2 dr.
\end{equation}
In addition, from \eqref{weightedhardy} we get
 \begin{equation}\label{eq2}
 I_2\geq \sum_{k=0}^{\infty} c_k g(N, k) \int_0^\infty r^{N-5} u_k^2 dr,
 \end{equation}
 where $g(N, k):=(N-4)^2/2+c_k+2(N-4)$. Since $\{c_k\}_{k\geq 0}$  is a nonnegative increasing sequence, it is easy to notice that the sequence $\{g(N, k)\}_{k\geq1}$ is positive and increasing for any $N\geq 3$. Therefore, we have
 $$g(N, k)\geq g(N, 1)=\frac{N^2-2N-2}{2}, \quad \forall k\geq 1.$$
 Since $c_0=0$ from \eqref{eq2} we obtain
\begin{equation}\label{eq3}
I_2\geq \frac{N^2-2N-2}{2} \sum_{k=0}^{\infty} c_k \int_0^\infty r^{N-5} u_k^2 dr
\end{equation}
 Summing up, from \eqref{eq1}, \eqref{eq3} and \eqref{gradform} we get
 \begin{equation}\label{eq4}
 \intr |\Delta u|^2 dx \geq \min\left\{\frac{N^2}{4}, \frac{N^2-2N-2}{2}\right\}\intr \frac{|\n u|^2}{|x|^2} dx.
 \end{equation}
 Since
 \begin{equation}\label{partres}
 \min\left\{\frac{N^2}{4}, \frac{N^2-2N-2}{2}\right\}=\left\{
 \begin{array}{ll}
   \frac{N^2}{4}, & N\geq 5 \\[5pt]
   3, & N=4 \\[5pt]
   \frac{1}{2} & N=3, \\ [5pt]
 \end{array}\right.
 \end{equation}
  inequality \eqref{optconst} is proven for any $N\geq 4$.

  For $N=3$ the final step of the argument above does not provide the optimal constant $C(3)$   since $1/2< C(3)=25/36$. In order to recover the constant $C(3)$ in the following we slightly modify the last part of the proof.

 First observe that the constant $N^2/4$ in \eqref{eq1} is optimal since the constant $(N-2)^2/4$ in inequality \eqref{WeightedRellich} is also optimal. This implies that
 $$C(N)\leq \frac{N^2}{4}, \quad \forall N\geq 3. $$
 and therefore, in view of \eqref{eq4} we obtain $C(N)=N^2/4$ for any $N\geq 5$. 
 
 For $N\in \{3, 4\}$ the minimum in \eqref{partres} is attained by $(N^2-2N-2)/2$ which is strictly smaller than $N^2/4$. In fact, due to this gap there is a coincidence that the minimum in \eqref{partres} for $N=4$ coincides with $C(4)$.

 In view of these considerations next we show how to recover the best constant $C(N)$ for $N\in \{3, 4\}$.   So, next we focus on  $N\in \{3, 4\}$.

 Observe that the term $\int_0^\infty r^{N-3} |u_k^{\prime}|^2 dr$ appears in both  $I_1$ and $I_2$. Next we want this term to be ``equally distributed" in $I_1$ and $I_2$ so that to contribute with the same constants in \eqref{eq1} and \eqref{eq3}.
  For that, first let $0<\eps< N^2/4$ which will be well precise later. Now we reconsider the terms $I_1$ and $I_2$ by  splitting  the right hand side of \eqref{finalLap} as $I_{1, \eps}+I_{2, \eps}$
  $$I_{1, \eps}:=\sum_{k=0}^{\infty} \left(\int_0^\infty r^{N-1}|u_k^{\prime\prime}|^2 dr + (N-1-\eps) \int_{0}^{\infty} r^{N-3}|u_k^{\prime}|^2 dr\right)$$
and
$$I_{2, \eps}:=\sum_{k=0}^{\infty} (2c_k+\eps) \int_{0}^{\infty} r^{N-3}|u_k^{\prime}|^2 dr+ \left(c_k^2+2c_k(N-4)\right)\int_0^\infty r^{N-5} u_k^2 dr.$$
 Again from \eqref{WeightedRellich} we obtain
 \begin{equation}\label{HL5}
 I_{1, \eps}:=\left(\frac{N^2}{4}-\eps\right) \sum_{k=0}^{\infty}  \int_{0}^{\infty} r^{N-3} |u_k^{\prime}|^2 dr.
 \end{equation}
 Applying \eqref{weightedhardy} and the fact that $c_0=0$ from the expression of $I_{2, \eps}$ we get
  \begin{equation}\label{HL6}
  I_{2, \eps} \geq  \sum_{k=1}^{\infty} c_k h(\eps, k) \int_0^\infty r^{N-5} u_k^2 dr,
  \end{equation}
 where $h(\eps, k):=(2 +\eps/c_k)(N-4)^2/4+c_k+2(N-4)$, for any $k\geq 1$. Since $c_k\geq N-1$ for any $k\geq 1$ we easily remark that the sequence $\{h(\eps, k)\}_{k\geq 1}$ is increasing. Therefore,
 $$h(\eps, k)\geq h(\eps, 1)=\left(2+\frac{\eps}{N-1}\right)\left(\frac{N-4}{2}\right)^2+3N-9, \quad \forall k\geq 1$$
 and it follows that
 \begin{equation}\label{HL7}
 I_{2, \eps}  \geq \left[\left(2+\frac{\eps}{N-1}\right)\left(\frac{N-4}{2}\right)^2+3N-9\right] \sum_{k=0}^{\infty} c_k
 \int_0^\infty r^{N-5} u_k^2 dr.
 \end{equation}
Next  we chose $\eps$ to obtain the same constant in both inequalities \eqref{HL5} and \eqref{HL7}, i.e.
$$ \frac{N^2}{4}-\eps=\left[\left(2+\frac{\eps}{N-1}\right)\left(\frac{N-4}{2}\right)^2+3N-9\right].$$
This is equivalent to
$$\eps(N)=\frac{(N-1)(-N^2+4N+4)}{N^2-4N+12}.$$
We then obtain
 \begin{equation}\label{HL8}
 \int_{\rr^N} |\Delta u|^2 dx \geq \left(\frac{N^2}{4}-\eps(N)\right)\int_{\rr^3}\frac{|\n \phi|^2}{|x|^2} dx.
 \end{equation}
 Since $\eps(4)=1$ and $\eps(3)=14/9$ we finally get the desired constants
 $$\frac{N^2}{4}-\eps(N)\Big|_{N=4}=3, \quad \frac{N^2}{4}-\eps(N)\Big|_{N=3}=\frac{25}{36}.$$
 We conclude that inequality \eqref{optconst} in Theorem \ref{HardyBilaplacian} holds also for $C(3)=25/36$ and $C(4)=3$.
\hfill $\square$

\begin{remark}
Notice also that the optimality of $C(N)=N^2/4$ for  $N\geq 5$ is hidden (and specified) in the proof of Theorem \ref{HardyBilaplacian} without the necessity of building a minimizing sequence.  \end{remark}

\section*{Proof of Theorem \ref{minimization}}
As we already mentioned in Remark \ref{rem1},  the proof of optimality is relevant only for $N\in \{3, 4\}$. However, for the sake of completeness, since our computations are slightly different than those in \cite{TertiZagra}, let us give a full dimensional proof.  
\paragraph{\bf Optimality (the cases $N\geq 5$).} Writing $u_\eps(x)=U_\eps(|x|) $, in view of \eqref{Bil}, since the spherical part is missing we obtain the simplified expression 
\begin{equation}\label{Bil2}
\int_{\rr^N}|\Delta u_\eps|^2 dx= |S^{N-1}| \left(\int_0^\infty r^{N-1}|U_{\eps}^{\prime\prime}(r)|^2 dr + (N-1) \int_{0}^{\infty} r^{N-3} |U_{\eps}^{\prime}(r)|^2 dr\right).
\end{equation}
and 
\begin{equation}\label{grad2}
\int_{\rr^N} \frac{|\nabla u_\eps|^2}{|x|^2} dx=|S^{N-1}|\int_{0}^\infty r^{N-3} |U_\eps^{\prime}(r)|^2 dr. 
\end{equation}
Then we have 
\begin{align}\label{id1}
\int_0^\infty  r^{N-3} |U_{\eps}^{\prime}(r)|^2 dr&=\left(-\left(\frac{N-4}{2}\right)+\eps\right)^2 \int_{0}^{\infty} r^{-1+2\eps } g^2(r)dr +\int_{0}^\infty r^{1+2\eps} g^{\prime}(r)^2 dr\nonumber\\
& + 2 \left(-\left(\frac{N-4}{2}\right)+\eps\right)\int_{0}^\infty r^{2\eps} g(r)g^{\prime}(r)dr\nonumber\\
&= \frac{1}{2\eps} \left(-\left(\frac{N-4}{2}\right)+\eps\right)^2 + \mathcal{O}(1). 
\end{align}
since $g^{\prime}$ is supported in the interval $[1, 2]$.  From the same reasons since 
$$U_\eps^{\prime\prime}(r)=\left(-\left(\frac{N-4}{2}\right)+\eps\right) \left(-\left(\frac{N-2}{2}\right)+\eps\right) r^{-N/2+\eps} g (r) + 
\chi_{[1, 2]} \mathcal{O}(1)$$
we obtain 
\begin{align}\label{id2}
\int_0^\infty r^{N-1}|U_{\eps}^{\prime\prime}(r)|^2 dr =\frac{1}{2 \eps} \left(-\left(\frac{N-4}{2}\right)+\eps\right)^2 \left(-\left(\frac{N-2}{2}\right)+\eps\right)^2 +\mathcal{O}(1). 
\end{align}
Due to \eqref{id1} and \eqref{id2} we successively obtain 
\begin{align*}
\frac{\int_{\rr^N}|\Delta u_\eps|^2 dx}{\int_{\rr^N}|\nabla u_\eps|^2/|x|^2 dx}&=\frac{ \left(-\left(\frac{N-4}{2}\right)+\eps\right)^2 \left(-\left(\frac{N-2}{2}\right)+\eps\right)^2+(N-1) \left(-\left(\frac{N-4}{2}\right)+\eps\right)^2+\mathcal{O}(\eps)}{\left(-\left(\frac{N-4}{2}\right)+\eps\right)^2 +\mathcal{O}(\eps)} \nonumber\\
& \searrow \frac{N^2}{4}=C(N), \quad \textrm{ as } \eps \searrow 0. 
\end{align*} 
The above limit also holds in the case $N=3$ but it does not provide the best constant $C(3)$. The case $N=4$ is not covered because of the nontermination $\frac{0}{0}$.  
\medskip
\paragraph{\bf Optimality (the cases $N\in \{3, 4\}$). }
As before we obtain 
\begin{align}\label{Bil2_lower}
\int_{\rr^N}|\Delta u_\eps|^2 dx&= \int_0^\infty r^{N-1}|U_{\eps}^{\prime\prime}(r)|^2 dr + (N-1+2c_1) \int_{0}^{\infty} r^{N-3} |U_{\eps}^{\prime}(r)|^2 dr\nonumber\\
&+ (c_1^2+2(N-4)c_1) \int_{0}^\infty r^{N-5} U_\eps^2(r) dr.  
\end{align}
and 
\begin{equation}\label{grad2_lower}
\int_{\rr^N} \frac{|\nabla u_\eps|^2}{|x|^2} dx=\int_{0}^\infty r^{N-3} |U_\eps^{\prime}(r)|^2 dr+c_1\int_{0}^\infty r^{N-5} U_\eps^2(r) dr. 
\end{equation}
Since 
$$\int_{0}^\infty r^{N-5} U_\eps^2(r)dr =\frac{1}{2\eps} +\mathcal{O}(1)$$
from \eqref{Bil2_lower} and \eqref{grad2_lower} we get that 
\begin{align*}
\frac{\int_{\rr^N}|\Delta u_\eps|^2 dx}{\int_{\rr^N}|\nabla u_\eps|^2/|x|^2 dx}&=\frac{ \left(\frac{N-4}{2}-\eps\right)^2 \left(\frac{N-2}{2}-\eps\right)^2+(N-1+2c_1) \left(\frac{N-4}{2}-\eps\right)^2+c_1^2+2(N-4)c_1 + \mathcal{O}(\eps)}{\left(-\left(\frac{N-4}{2}\right)+\eps\right)^2+ c_1+\mathcal{O}(\eps)} \nonumber\\
& \searrow\frac{ \left(\frac{N-4}{2}\right)^2 \left(\frac{N-2}{2}\right)^2+(N-1+2c_1) \left(\frac{N-4}{2}\right)^2+c_1^2+2(N-4)c_1}{\left(\frac{N-4}{2}\right)^2+ c_1}, \quad \textrm{ as } \eps \searrow 0. 
\end{align*} 
Since 
$$\frac{ \left(\frac{N-4}{2}\right)^2 \left(\frac{N-2}{2}\right)^2+(N-1+2c_1) \left(\frac{N-4}{2}\right)^2+c_1^2+2(N-4)c_1}{\left(\frac{N-4}{2}\right)^2+ c_1}=\left\{
\begin{array}{cc}
3, & \textrm{ if } N=4, \\[3pt] 
\frac{25}{36}, &  \textrm{ if } N=3. \\  
\end{array}\right.$$
The proof of optimality is proved. 

\medskip
\paragraph{\bf The non-attainability of the best constant $C(N)$, $N\geq 3$.}

	The non-attainability follows the lines of the proof of Theorem \ref{HardyBilaplacian}.  Indeed, assuming that $C(N)$ is attained then it is necessary to have equality in inequalities \eqref{WeightedRellich}-\eqref{weightedhardy} for any $u_k$ in the decomposition of $u$. Remark that inequality \eqref{WeightedRellich} is also a consequence of the identity 
	
	\begin{equation}\label{identi}
	\int_0^\infty r^{N-1} |u_k^{\prime \prime}|^2 dr - \frac{(N-2)^2}{4} \int_0^\infty r^{N-3} |u_k^{\prime}|^2 dr= \int_0^\infty \left|\left(r^{\frac{N-2}{2}}u_k^{\prime}\right)^{\prime}\right|^2 r dr.
	\end{equation}
In view of \eqref{identi} we obtain  that equality in  \eqref{WeightedRellich}	 is achieved  if
$$\left(r^{\frac{N-2}{2}}u_k^{\prime}\right)^{\prime}=0$$
which leads to the family of solutions
$$u_k=m_k r^{-\frac{N-4}{2}}+ n_k,$$
	for some real constants $m_k, n_k$, with the fundamental system of solutions given by $\{r^{-\frac{N-4}{2}}, 1\}$.  
	Observe that $u_k=1$ is not possible since constant functions are not admissible for inequality \eqref{weightedhardy}. 
	On the other hand, $u_k=r^{-\frac{N-4}{2}}$ is not admissible either because none of the terms in \eqref{WeightedRellich} is integrable. In consequence   the constant $C(N)$ is not attained. \hfill $\square$



\subsection*{Acknowledgements}

This work was partially supported by  a grant of Ministery of Research and Innovation, CNCS-UEFISCDI, project number PN-III-P1-1.1-TE-2016-2233, within PNCDI III and by 
  a Young Researchers Grant awarded by The Research Institute of the University of Bucharest (ICUB).

\bibliographystyle{amsplain}
\bibliography{biblio1}
\end{document}